\documentclass[12pt]{article}
\usepackage[utf8]{inputenc}
\usepackage[english]{babel}
\usepackage{amsthm,amssymb, amsmath, amsfonts, amscd}
\usepackage{graphicx}
\usepackage{textcomp}
\usepackage{euscript}
\usepackage{multirow}
\usepackage[matrix , arrow, curve]{xy}
\def\R{\mathbb{R}}
\def\beginproof{{\bf Proof. }}
\def\endproof{\hfill$\square$\medskip}
\newtheorem{thm}{Theorem}
\newtheorem{lemma}{Lemma}
\newtheorem{proposition}{Proposition}
\newtheorem{definition}{Definition}

\theoremstyle{remark}

\begin{document}
\centerline{\large \bf Multiplicity of positive solutions for}
\centerline{\large \bf a critical quasilinear Neumann problem}
\medskip
\centerline{Aleksandr Enin\footnote{Saint Petersburg State University, St.~Petersburg, Russia; Saint Petersburg Electrotechnical University "LETI", St.~Petersburg, Russia; Aleksandr.Enin@outlook.com}}

\begin{abstract}
We establish the multiplicity of positive solutions to a quasilinear Neumann problem in expanding balls and hemispheres with critical exponent in the boundary condition.
\end{abstract}
\section{Introduction}
We consider the following problem
\begin{equation}\label{expanding_problem}
\left\{
\begin{array}{lll}
&\Delta_p u := \text{div}(|\nabla u|^{p-2} \nabla u) = |u|^{p-2}u \quad \mbox{in } B_R, \\\\
&|\nabla u|^{p-2} \langle\nabla u; {\bf n} \rangle = |u|^{q-2}u \quad \mbox{on } S_R,\\\\
& u>0 \quad \mbox{in } B_R,\\
\end{array}
\right.
\end{equation}
where $B_R$ and $S_R$ are the ball and the sphere with radius $R$ respectively in $\mathbb{R}^n$. Here $1 < p < n$  and $q = p^{**} = \frac {(n-1)p} {(n - p)}$ is the critical exponent for the trace embedding.

We establish the multiplicity effect for weak solutions to \eqref{expanding_problem}. Namely we prove that the number of positive rotationally non-equivalent solutions is unbounded as $R \to \infty$.

The effect of multiplicity was discovered by Coffman \cite{Coffman} who considered the Dirichlet problem
\begin{equation}
\label{dirichlet_problem}
\left\{
\begin{array}{lll}
&-\Delta_p u = |u|^{q-2}u \quad \mbox{in } \Omega_R, \\\\
& u=0 \quad \mbox{on } \partial\Omega_R,\\\\
& u>0 \quad \mbox{in } \Omega_R, \\
\end{array}
\right.
\end{equation}
where $\Omega_R$ is the annulus $B_R \setminus B_{R-1} \subset \R^n$ for $n=2$ and $p=2$.  The problems \eqref{expanding_problem} and \eqref{dirichlet_problem} were studied later by many authors for subcritical $q$ (see e.g. \cite{Scheglova, Nazarov&Enin, Kolonitskii, Li, Byeon}).
In \cite{Wang} the multiplicity result was obtained for the Neumann problem
\begin{equation*}
\left\{
\begin{array}{lll}
&-\Delta u + \lambda u = |u|^{p^*-2}u \text{ in } \Omega,\; \\\\
&\frac{\partial u}{\partial \nu} = 0 \text{ on } \partial\Omega, \\\\
& u > 0 \quad \mbox{in } \Omega, \\
\end{array}
\right.
\end{equation*}
where $\Omega$ satisfies some symmetry conditions and $p^*$ is the critical exponent for the Sobolev trace embedding.

One can easily show that after suitable rescaling solutions of \eqref{expanding_problem} are solutions to the following problem:
\begin{equation}\label{problem}
\left\{
\begin{array}{lll}
&\Delta_p u = \lambda |u|^{p-2}u \quad \mbox{in } B, \\\\
&|\nabla u|^{p-2} \langle\nabla u; {\bf n} \rangle = |u|^{q-2}u \quad \mbox{on } S,\\\\
& u>0 \quad \mbox{in } B,\\
\end{array}
\right.
\end{equation}
where $B=B_1$, $S=S_1$ and $\lambda(R) = R^{p}$ as $R\to\infty$. 

We look for distinct solutions of the problem \eqref{problem} by minimizing the functional
\begin{equation}
\label{energy_functional}
I^\lambda[u] := \frac{\|\nabla u\|^p_{L_p(B)} + \lambda \|u\|^p_{L_p(B)}}{\|u\|^p_{L_q(S)}}
\end{equation}
on different subsets of $W^1_p(B)$.

In order to construct solutions to problem \eqref{problem} let us introduce the following notation:
\begin{definition}
Let $A \subset S$ and $\varkappa>0$. We denote by $A^\delta$ $\varkappa$-neighborhood of a set $A$, i.e.
\begin{equation*}
A^\varkappa = \{z \in S \; | \; \text{dist}(z,A)\leq \delta\}.
\end{equation*}
\end{definition}
The following definition was introduced in \cite{Byeon}:
\begin{definition}
Let $G$ be a closed subgroup of $O(n)$. We call set $A \subset S$ a locally minimal orbital set under the action of $G$ if $A$ is invariant under the action of $G$ and satisfies the following conditions:
\begin{itemize}
\item for any $x \in A$ the orbit $Gx$ is a discrete set and $m(A):=|Gx|$ is independent of $x$.
\item there exists $\varkappa > 0$ such that for any $y \in A^\varkappa \setminus A$ and $x\in A$, we have $|Gx| < |Gy|$.
\end{itemize}

\end{definition}

We denote as $m(G)$ the number of elements in the minimal orbit of $G$ and $K(n,p)$ stands for the best Sobolev trace constant in half-space defined as
\begin{equation*}
K(n,p) = \inf_{v \in C_c^\infty(\overline{\mathbb{R}^n_+})\backslash \{0\}} \frac{ \|\nabla v\|^p_{L_p(\mathbb{R}^n_+)}}{\|v(\cdot, 0)\|^p_{L_q(\mathbb{R}^{n-1})}}.
\end{equation*}
The value of $K(n,p)$ is calculated explicitly in \cite{Escobar} for $p=2$ and \cite{Nazaret_embedding} for arbitrary $p$.

We consider local minimizers of functional \eqref{energy_functional} on sets
\begin{equation*}
X_G(A, \beta) = \{u \in W^1_p(B) \; | \; u(gx) \equiv u(x) \; \forall g\in G, \; \|u\|_{L_q(S)} = 1, \; \|u\|^q_{L_q(A^{\varkappa})} \geq 1 - \beta\},
\end{equation*}
where $G$ is some closed subgroup of $O(n)$, $A$ is a locally minimal orbital set and $\beta$ is some small parameter that we will choose later. We denote $X_G(A,\beta)$ by $X$ if it does not lead to confusion.

The structure of the paper is as follows. In Section 2 we prove some auxiliary lemmas and in Section 3 we establish main multiplicity results.
\section{Auxiliary lemmas}

The following fact is well known and will be given here without a proof.
\begin{proposition}
\label{Gateaux differential}
The functional $I^\lambda[u]$ is Gateaux differentiable and for any $h\in W^1_p(B)$
\begin{equation*}
\begin{aligned}
DI^\lambda[u](h) &= p \int_B |\nabla u|^{p-2} \nabla u \cdot \nabla h \, dx \frac{1}{\|u\|_{L_q(S)}^p} - p \int_B |\nabla u|^p dx \int_{S} |u|^{q-2} uh \, dS \frac{1}{\|u\|_{L_q(S)}^{p+q}} \\
&- p\lambda\int_B |u|^p dx \int_{S} |u|^{q-2} uh \, dS \frac{1}{\|u\|_{L_q(S)}^{p+q}} + p \lambda \int_B |u|^{p-2} uh \, dx  \frac{1}{\|u\|_{L_q(S)}^p}.
\end{aligned}
\end{equation*}
\end{proposition}

\begin{lemma}
\label{delta_estimate}
Let $u^\lambda_j \in W^1_p(B)$ be a bounded Palais-Smale sequence for $I^\lambda$ at the level $c > 0$. Then there is $u^\lambda_0 \in W^1_p(B)$ such that up to subsequence $u^\lambda_j \rightharpoondown u^\lambda_0$ and
\begin{eqnarray}
\label{mu_estimate}
|\nabla u_j^\lambda|^p dx \rightharpoondown \mu \geq |\nabla u^\lambda_0|^p dx + \sum_k \mu_k \delta (x-x_k), \\
\label{PS_estimate}
|u_j^\lambda|^q dS \rightharpoondown \nu = |u^\lambda_0|^q dS + \sum_k \nu_k \delta (x-x_k),
\end{eqnarray}
where $\delta (x-x_k)$ are delta measures at some points $x_k$ in $S$ and $\mu_k \geq K(n,p)\nu_k^{\frac{p}{q}}$.
Furthermore, either $\nu_k = 0$ or $\nu_k \geq (c^{-1}\cdot K(n,p))^{\frac{q}{q-p}}\nu(S)$.
\end{lemma}
\beginproof
Since $\{u^\lambda_j\}$ is bounded in $W^1_p(B)$, the relations \eqref{mu_estimate} and \eqref{PS_estimate} follow by the Lions concentration-compactness principle \cite{Lions}. Since $I^\lambda$ is homogeneous we can assume without loss of generality that $\|u_j^\lambda\|_{L_q(S)} = 1$ and $\nu(S) = 1$.
Next we use the argument from \cite{Nazaret, Bonder}:
Let us fix $x_k$ from \eqref{mu_estimate} and \eqref{PS_estimate}. We choose $\varphi \in C_c^\infty(\R^n)$ such that
\begin{equation*}
\varphi = 1 \text{ in } B(x_k, \varepsilon), \quad \quad \varphi = 0 \text{ in } \mathbb{R}^n \setminus B(x_k, 2\varepsilon), \quad |\nabla \varphi| \leq \frac{C}{\varepsilon}.
\end{equation*}
Since $DI^\lambda[u^\lambda_j] \to 0$ we obtain
\begin{equation*}
\lim\limits_{j \to \infty} <DI^\lambda[u^\lambda_j], \varphi u^\lambda_j> = 0
\end{equation*}
Then
\begin{equation}
\label{atom_estimate}
\lim_{j \to \infty} \int_B |\nabla u^\lambda_j|^{p-2} \nabla u^\lambda_j \cdot \nabla \varphi_j u^\lambda_j\, dx = c \int_S \varphi d\nu - \int_B  \varphi d\mu - \lambda \int_B |u^\lambda_0|^p \varphi dx.
\end{equation} 

One can estimate the left hand side as follows:
\begin{equation*}
\begin{aligned}
0 &\leq \left| \lim_{j \to \infty} \int_B |\nabla u^\lambda_j|^{p-2} \nabla u^\lambda_j \cdot \nabla \varphi u^\lambda_j\, dx \right| \\
&\leq \lim_{j \to \infty} \left(\int_B |\nabla u^\lambda_j|^p dx\right)^\frac{p-1}{p} \left(\int_B |\nabla \varphi|^p |u^\lambda_j|^p dx \right)^\frac{1}{p}	\\
&\leq C \left(\int_{B(x_k, 2\varepsilon)}|\nabla \varphi|^p |u^\lambda_0|^p dx \right)^\frac{1}{p}\\
&\leq C \left(\int_{B(x_k, 2\varepsilon)} |\nabla \varphi|^n dx\right)^\frac{1}{n} \left(\int_{B(x_k, 2\varepsilon)} |u^\lambda_0|^\frac{np}{n-p}dx\right)^\frac{n-p}{pn} \\
&\leq C \left(\int_{B(x_k, 2\varepsilon)} |u^\lambda_0|^\frac{np}{n-p}dx\right)^\frac{n-p}{pn} \to 0 \text{ as } \varepsilon \to 0.
\end{aligned}
\end{equation*}

Taking the limit in \eqref{atom_estimate} we get $\nu_k = c^{-1} \mu_k \geq c^{-1} K(n,p) \nu_k^\frac{p}{q}$. This means either $\nu_k \geq (c^{-1} K(n,p))^{\frac{q}{q-p}}$ or $\nu_k = 0$.
\endproof

\begin{lemma}
\label{concentration_lemma}
Let $u^\lambda \in X$ be a sequence such that $I^\lambda[u^\lambda] \leq K(n,p) m(A)^{1-\frac{p}{q}}$.
Then there is a $\beta_0 > 0$ such that for any $\beta \leq \beta_0$ there is $x_0 \in S$ such that we have up to subsequence the following weak convergence in the sense of measures as $\lambda \to \infty$:
\begin{equation} 
\label{weak_convergence}
|u^\lambda|^q dS \rightharpoondown \sum_{x_k \in Gx_0} \frac{1}{m(A)} \delta(x-x_k).
\end{equation} 
\end{lemma}
\beginproof
Since $\|u^\lambda\|^p_{W^1_p} \leq I^\lambda[u^\lambda] \leq K(n,p) m(A)^{1-\frac{p}{q}}$ by the Lions concentration-compactness principle we get 
\begin{eqnarray*}
|\nabla u^\lambda|^p dx \rightharpoondown \mu \geq |\nabla u_0|^p dx + K(n,p) \sum_k \nu_k^{\frac{p}{q}} \delta (x-x_k), \\
|u^\lambda|^q dS \rightharpoondown \nu = |u_0|^q dS + \sum_k \nu_k \delta (x-x_k),
\end{eqnarray*}
where $\delta (x-x_k)$ are delta measures at some points $x_k$ in $S$.

Since $\lambda \|u^\lambda\|^p_{L_p(B)}$ is uniformly bounded, we have $u^\lambda \to 0$ in $L_p(B)$  so $u_0 = 0$. Combining the above with the fact that $u^\lambda$ are invariant with respect to $G$ we get:
\begin{align}
\label{J_estimate}
\lim_{\lambda \to \infty} I^\lambda[u^\lambda] = \mu(B) \geq K(n,p) \sum_k \nu_k^\frac{p}{q} &= K(n,p) \sum_j |Gx_j| \nonumber
 \left(\frac{\tilde{\nu_j}}{|Gx_j|}\right)^\frac{p}{q} \\ &= K(n,p) \sum_j |Gx_j|^{1 - \frac{p}{q}} \tilde{\nu_j}^\frac{p}{q}. 
\end{align}
Here $j$ goes over different classes of equivalence of $x_k$, and $\tilde{\nu}_j = |Gx_j| \nu_j$ is a total contribution of that class to $\nu(\partial\Omega)$. The second equality is due to the fact that $u^\lambda$ are $G$-invariant, so for every $x_k$ there are $|Gx_k|$ $\delta$-functions with the same coefficient.

Since $p<q$ we have $a^{\frac{p}{q}} + b^{\frac{p}{q}} > (a+b)^\frac{p}{q}$, for any $a > 0,b > 0$.  Recalling that $A$ is a locally minimal orbital set  we can write
\begin{equation}
\label{concentration_estimate}
\begin{aligned}
\mu(B) &\geq K(n,p) m(A)^{1 - \frac{p}{q}} \sum_{j \;: \;x_j \in A} \tilde{\nu}_j^{\frac{p}{q}} + K(n,p) \sum_{i \;:\; x_i \not\in A} |Gx_i|^{1-\frac{p}{q}} \tilde{\nu}_i^{\frac{p}{q}} \geq \\
&\geq K(n,p) (m(A)^{1-\frac{p}{q}}\alpha^{\frac{p}{q}} + m(G)^{1 - \frac{p}{q}} (1-\alpha)^{\frac{p}{q}}),
\end{aligned}
\end{equation}
where $1  -\beta \leq \alpha \leq 1$ (we recall that $m(G)$ is the number of elements in the minimal orbit of $G$).

It's easy to see that the right hand side of \eqref{concentration_estimate} is a concave function of $\alpha$. That means that if $\beta$ is small enough, then the right hand side is a decreasing function, which achieves it's minimum of $K(n,p)m(A)^{1-\frac{p}{q}}$ at $\alpha = 1$.

Since by assumption $\mu(B) = \lim_{\lambda \to \infty} I^\lambda[u^\lambda] \leq K(n,p)m(A)^{1-\frac{p}{q}}$ we conclude that $\alpha = 1$. Recalling that for $u\in X$ $\|u\|_{L_q(S)} = 1$ we get \eqref{weak_convergence}.
\endproof

From now on we always assume that $\lambda$ is fixed and whenever there is a limit it is taken over $j\to\infty$ unless specified otherwise.

\begin{lemma}
\label{minimum_attained}

The minimum of $I^\lambda$ on $X$ is attained if $\lambda$ is large enough and
\begin{equation*}
\inf\limits_{u\in X} I^\lambda\left[u\right] < K(n,p)m(A)^{1-\frac{p}{q}}.
\end{equation*}
\end{lemma}
\beginproof
The Ekeland’s variational principle \cite{Ekeland} provides the existence of a minimizing sequence $u^\lambda_j \in X$ such that $I'[u^\lambda_j] \to 0$. Since $u^\lambda_j$ is a Palais-Smale sequence at the level $\inf\limits_{u\in X} I^\lambda[u] < K(n,p)m(A)^{1-\frac{p}{q}}$, Lemma \ref{delta_estimate} gives the estimate on any non-zero $\nu_k$ in \eqref{PS_estimate}:
\begin{equation}
\label{nu_estimate}
\nu_k > m(A)^{-\left(1 - \frac{p}{q}\right) \frac{q}{q-p}} = \frac{1}{m(A)}.
\end{equation}

Suppose that there is a $\delta$-function outside of $A$. 
From \eqref{weak_convergence} follows that for large $\lambda$ almost all of $\nu(S)$ mass is concentrated in a $\varkappa$-neighbourhood of $A$, and according to \eqref{nu_estimate} there are no $\delta$-functions outside of that neighbourhood.

Let us suppose that there is a $\delta$-function at $x_k\in A^\varkappa$ with weight $\nu_k$. Since $A$ is a locally minimal orbital set, we know that $|Gx_k| \geq m(A)$. Now from \eqref{J_estimate} and \eqref{nu_estimate} we derive
\begin{equation}
\lim\limits_{j \to \infty} I^\lambda[u^\lambda_j] \geq K(n,p) |Gx_k|\left(\frac{1}{m(A)}\right)^{\frac{p}{q}} = K(n,p) m(A)^{1 - \frac{p}{q}},
\end{equation}
which is a contradiction.

From that follows that for $u^\lambda_0$ in \eqref{weak_convergence} $\|u^\lambda_0\|_{L_q(S)} = \|u^\lambda_j\|_{L_q(S)} = 1$. It is well known, that weak convergence and convergence of norms implies strong convergence in uniformly convex Banach space (e.g. \cite[Proposition~3.32]{Brezis}), and that completes our proof. That way $u^\lambda_0 \in X$ and $I^\lambda[u^\lambda_0]$ attains minimal value. 
\endproof

\section{Main results}

\begin{lemma}
\label{upper_bound}
Let $G = H\times O(n-k)$ where $H$ is a finite subgroup of $O(k)$ and $A \subset \mathbb{R}^k$ is a minimal orbital set under the action of $H$.

Then for any fixed $\beta$, $\lambda$ large enough and $p\leq \frac{n+1}{2}$ we have
\begin{equation}
\label{norm_estimate}
\inf\limits_{u \in X} I^\lambda[u^\lambda] < K(n,p) m(A)^{1 - \frac{p}{q}}.
\end{equation}
\end{lemma}
\beginproof
Let $x_0 \in Gx_0$ be a point in $A\times\{0\}$.
As was shown in \cite{Nazarov_embedding} (see also \cite{Bonder_embedding}) there is a function $u_R$ in $W^1_p(B_R)$ supported in a small ball around $Rx_0$ and axially symmetric with respect to the axis $Ox_0$, such that $\|u_R\|^p_{W^1_p(B_R)} < K(n,p) \|u\|_{L_q(S_R)}^p$.

Now we construct the function
\begin{equation*}
v_R(x) = \sum_{g \in H} u_R(gx).
\end{equation*}

It is easy to see that $v_R$ is $G$-invariant and
\begin{equation*}
\frac{\|v_R\|^p_{W^1_p(B_R)}}{\|v_R\|^p_{L_q(S_R)}} = m(A)^{1-\frac{p}{q}} \frac{\|u_R\|^p_{W^1_p(B_R)}}{\|u_R\|^p_{L_q(S_R)}} < K(n,p) m(A)^{1-\frac{p}{q}}.
\end{equation*}
By rescaling we obtain \eqref{norm_estimate}.
\endproof

\begin{thm}
\label{G-invariant_solution}
Let $p \leq \frac{n+1}{2}$ and let $G$ be as in Lemma \ref{upper_bound}. Suppose that $A \subset \mathbb{R}^k$ is some locally minimal orbital set of H. Then there is $\lambda_0 > 0$ such that for any $\lambda > \lambda_0$ there is a $G$-invariant solution of problem \eqref{problem} such that it concentrates at $|Gx_0|$ points in the $Gx_0$ for some $x_0 \in A\times \{0\}$, i.e.
\begin{equation*}
\frac{|u^\lambda|^q}{\|u^\lambda\|_{L_q(S)}} \rightharpoondown \sum_{k=1}^{|G(x_0)|} \frac{1}{|G(x_0)|}\delta(x-x_k) \quad\quad \text{ as } \quad \lambda \to\infty.
\end{equation*}
\end{thm}
\beginproof
According to Lemmas \ref{upper_bound} and \ref{minimum_attained} there is a minimizer $u \in X$ such that it is concentrated around $m(A)$ points of $A\times\{0\}$.
Lemma \ref{concentration_lemma} implies that if $\lambda$ is large enough the constraint $\|u\|^q_{L_q(A^{\delta})} > 1 - \beta$ is non-active and does not produce a Lagrange multiplier. Since $I^\lambda[u] = I^\lambda[|u|]$ we can assume that $u$ is non-negative. Since $u$ is a local minimizer, we get for $\mu = I^\lambda[u]$ (see Proposition \ref{Gateaux differential}):
\begin{equation*}
\int_B |\nabla u|^{p-2} \nabla u \cdot \nabla h \, dx + \lambda \int_B |u|^{p-2} uh \, dx - \mu \int_{S} |u|^{q-2} uh \, dS = 0 \quad\quad \forall h \in L_G,
\end{equation*}
where 
\begin{equation*}
L_G = \{h\in W^1_p(B)\;|\; h(gx) = h(x) \; \forall g\in G\}.
\end{equation*}
Due to the principle of symmetric criticality \cite{Symmetric_criticality} $u$ is a solution to the problem 
\begin{equation*}
\left\{
\begin{array}{ll}
&\Delta_p u := \lambda |u|^{p-2}u \quad \mbox{in } B, \\\\
&|\nabla u|^{p-2} \langle\nabla u; {\bf n} \rangle = \mu |u|^{q-2}u \quad \mbox{on } S,
\end{array}
\right.
\end{equation*}
Since $u\geq 0$ in $B$ we can apply the Harnack inequality (see \cite{Trudinger}, \cite{Serrin}) and get the positivity of our solution. Since the boundary condition is not homogeneous, it's easy to show that $\mu^{\frac{1}{p-q}} u$ is a solution for problem \eqref{problem}.
\endproof

\begin{thm}
For any $N>0$ there is $\lambda_0 > 0$ such that for every $\lambda > \lambda_0$ problem \eqref{problem} has at at least $N$ distinct solutions.
\end{thm}
\beginproof
Let us look at the following decomposition of $\mathbb{R}^n$:
\begin{equation*}
\mathbb{R}^n = \left(\mathbb{R}^2\right)^l \times \mathbb{R}^m.
\end{equation*}
Here $l\geq 1$, $m \geq 0$. We denote variables in $\mathbb{R}^n$ by $x$, in $\mathbb{R}^2$ by $y$ and in $\mathbb{R}^m$ by $z$. This way,
\begin{equation*}
x = (y_1, y_2, \ldots, y_l; z).
\end{equation*}

We introduce the group $G_{k,l} = H_{k,l} \times O(m)$ where $H_{k,l}$ is generated by rotations of every $y_i$ by $\frac{2\pi}{k}$ and by transpositions of $y_i$ and $y_j$ for every $i$ and $j$. 

Let $A$ be a globally minimal orbital set for the action of $H_{k,l}$. One can easily check that $A \times \{0\}$ is a locally minimal orbital set for $G_{k,l}$.



Now we show that for $l \geq 1$ and $k>2$ the minimizers will be non-equivalent. In order to do that we analyse minimal orbits of $H_{k,l}$. The simple calculation yields that a minimal orbit would be of a point $(y,0,\ldots, 0) \in \R^{2l}$ where $y \in \R^2$ and it consists of $k\cdot l$ points. Knowing the structure of the minimal orbits we can deduce that minimizers would be different for different pairs of $(k,l)$ and $(k',l')$.
\endproof

Now we consider an analogue of the problem \eqref{problem} in an $n$-dimensional hemisphere.


To prove the multiplicity result we only need to modify lemma \ref{upper_bound} by using the existence result from \cite{NR}.
\begin{lemma}
Let $n\geq 5$ and let $B$ be an $n$-dimensional hemisphere. 
Let $G = H\times O(n-k)$ where $H$ is a finite subgroup of $O(k)$ such that $A$ is a minimal orbital set under the action of $H\times \{0\}$.

Then for any fixed $\beta$, $\lambda$ large enough and $2 < p \leq \frac{n+2}{3}$ we have
\begin{equation*}
\inf\limits_{u \in X} I^\lambda[u^\lambda] < K(n,p) m(A)^{1 - \frac{p}{q}}.
\end{equation*}
\end{lemma}

Repeating the previous arguments we get the following theorem:

\begin{thm}
Let $n \geq 5$, and let $B$ be an $n$-dimensional hemisphere, $2 < p \leq \frac{n+2}{3}$. Then for any $N > 0$ there is a $\lambda_0 > 0$ such that for any $\lambda > \lambda_0$ problem \eqref{problem} has at least $N$ rotationally non-equivalent solutions.
\end{thm}

\subsection*{Acknowledgements} 
I am grateful to A.I. Nazarov for stimulating discussions.

The work has been supported by RFBR №14-01-00534.

A part of this paper was written during the visit of author to University of Cologne, supported by joint program of DAAD	and St. Petersburg State University "Dmitrij Mendeleev". I am grateful to professor Bernd Kawohl for the hospitality during this visit.

\end{document}